\documentclass[preprint,10pt]{elsarticle}

\usepackage{tikz}
\colorlet{myblue}{blue!60!black}
\definecolor{bleudefrance}{rgb}{0.19, 0.55, 0.91}

\usepackage{hyperref}
\usepackage{amsmath,amssymb,amsthm}
\usepackage{dsfont}  

\newcommand{\C}{\mathbb C}
\newcommand{\R}{\mathbb R}
\newcommand{\Q}{\mathbb Q}
\newcommand{\Z}{\mathbb{Z}}
\newcommand{\N}{\mathbb{N}}

\def\boldb{ \boldsymbol{b}}
\def\bc{ \boldsymbol{c}}

\def\bv{ \boldsymbol{v}}

\newcommand{\kernel}{\mathrm{ker}} 

\newcommand{\sgn}{\mathrm{sgn}} 

\makeatletter
\def\moverlay{\mathpalette\mov@rlay}
\def\mov@rlay#1#2{\leavevmode\vtop{%
   \baselineskip\z@skip \lineskiplimit-\maxdimen
   \ialign{\hfil$\m@th#1##$\hfil\cr#2\crcr}}}
\newcommand{\charfusion}[3][\mathord]{
    #1{\ifx#1\mathop\vphantom{#2}\fi
        \mathpalette\mov@rlay{#2\cr#3}
      }
    \ifx#1\mathop\expandafter\displaylimits\fi}
\makeatother

\newtheorem*{thm-nonumber}{Theorem}			
\newtheorem{theorem}{Theorem}
\newtheorem{corollary}[theorem]{Corollary}
\newtheorem{lemma}[theorem]{Lemma}
\newtheorem*{prop-nonumber}{Proposition}			
\newtheorem{proposition}[theorem]{Proposition}
\newtheorem{example}[theorem]{Example}
\newtheorem*{remark}{Remark}		
\newtheorem{Conjecture}{Conjecture}

\usepackage{color}

\begin{document}

\begin{frontmatter}



\title{Bases of complex exponentials with restricted supports}


\author[MIDS]{Dae Gwan Lee\corref{cor1}}
\ead{daegwans@gmail.com}
\cortext[cor1]{Corresponding author}

\author[MIDS]{G\"otz E.~Pfander}
\ead{pfander@ku.de}

\author[GMU]{David Walnut}
\ead{dwalnut@gmu.edu}


\address[MIDS]{Mathematical Institute for Machine Learning and Data Science (MIDS), Katholische Universit\"at Eichst\"att--Ingolstadt, Goldknopfgasse 7, 85049 Ingolstadt, Germany}

\address[GMU]{Department of Mathematical Sciences, George Mason University, Fairfax, VA 22030, USA}




\begin{abstract}
The complex exponentials with integer frequencies form a basis for the space of square integrable functions on the unit interval. We analyze whether the basis property is maintained if the support of the complex exponentials is restricted to possibly overlapping subsets of the unit interval.  We show, for example, that if $S_1, \ldots, S_K \subset [0,1]$ are finite unions of intervals with rational endpoints that cover the unit interval, then there exists a partition of $\Z$ into sets $\Lambda_1, \ldots, \Lambda_K$ such that $\bigcup_{k=1}^K \{ e^{2\pi i \lambda (\cdot)} \chi_{S_k} : \lambda \in \Lambda_k \}$ is a Riesz basis for $L^2[0,1]$. Here, $\chi_S$ denotes the characteristic function of $S$.
\end{abstract}



\begin{keyword}
complex exponentials \sep
Riesz bases \sep
support restriction \sep
spectrum 


\MSC 42C15
\end{keyword}

\end{frontmatter}



\section{Introduction and Main Results}
\label{sec:intro}

For a measurable set $S \subset \R$ and a discrete set $\Lambda \subset \R$, we define $\mathcal{E} (S,\Lambda) := \{ e^{2\pi i \lambda (\cdot)} \chi_{S} : \lambda \in \Lambda \}$ which is the set of all exponential functions with frequencies in $\Lambda$ considered on the domain $S$.  
A {\it Riesz basis} for a separable Hilbert space $\mathcal{H}$ is a sequence of the form $\{ U e_n \}_{n \in \Z}$, where $\{ e_n \}_{n \in \Z}$ is an orthonormal basis for $\mathcal{H}$ and $U : \mathcal{H} \rightarrow \mathcal{H}$ is a bijective bounded operator.
Equivalently, Riesz bases are characterized as complete Riesz sequences\footnote{A sequence $\{ f_n \}_{n\in\Z} \subset \mathcal{H}$ is called a {\it Riesz sequence} if there exist constants $0 < A \leq B < \infty$ such that
$A \, \| c \|_{\ell_2}^2
\leq \| \sum_{n\in\Z} c_n f_n \|^2
\leq B \, \| c \|_{\ell_2}^2$
for all $c=\{c_n\}_{n\in\Z} \in \ell_2 (\mathbb Z)$.} (see, e.g., \cite[Chapter 3.6]{Ch16}).

Most of the results on exponential bases deal with exponential functions that are defined on the full domain of the involved space. For instance, if one speaks of exponential functions in the space $L^2(S)$, then functions of the form $t \mapsto e^{2\pi i \lambda t}$ with $\lambda \in \R$, defined on the full domain $S$, are usually considered.

In this paper, we consider the situation where exponential functions in $L^2(S)$ are restricted to some different subsets of $S$.
More precisely, we are interested in conditions on $S_1, \ldots, S_K \subset \R$ and $\Lambda_1, \ldots, \Lambda_K \subset \R$ which are necessary/sufficient for $\bigcup_{k=1}^K \mathcal{E}(S_k, \Lambda_k)$ to be a Riesz basis for $L^2( \bigcup_{k=1}^K S_k )$.

As our first main result, we prove that if each $S_k \subset [0,1)$ is a finite union of intervals with rational endpoints, then one can find pairwise disjoint sets $\Lambda_1, \ldots, \Lambda_K \subset \Z$ such that $\bigcup_{k=1}^K \mathcal{E}(S_k, \Lambda_k)$ is a Riesz basis for $L^2( \bigcup_{k=1}^K S_k )$.

\begin{theorem}\label{thm:rational-endpoints}
Let $K, N \in \N$ with $K \leq N$.
Let $I_1, \ldots, I_K$ be distinct intervals from $[\frac{\ell}{N}, \frac{\ell+1}{N})$, $\ell = 0, \ldots, N{-}1$, 
and for each $k$, let $S_k$ be a union of subcollection of $I_1, \ldots, I_K$, including $I_k$.
Then there exists a permutation\footnote{We denote by $\boldsymbol{\mathrm{S}}_K = \boldsymbol{\mathrm{S}} (\{1 , \ldots, K \} )$ the symmetric group of $\{1, \ldots, K \}$, which consists of all the $K!$ permutations of $1, \ldots, K$.} $\rho \in \boldsymbol{\mathrm{S}}_K$ such that
\[
\bigcup_{k=1}^{K} \mathcal{E} ( S_k,  N\Z {+} \rho(k) )
=\bigcup_{k=1}^{K}\big\{ e^{2\pi i \lambda\, (\cdot)}\chi_{S_k} : \lambda\in N\Z {+} \rho(k) \big\}
\]
is a Riesz basis for $L^2(\bigcup_{k=1}^{K} S_k)$.
\end{theorem}

We point out that if, for instance, $S_1 = [0,1)$, $S_2 = \cdots = S_N = [0,\frac{1}{N})$, then necessarily the system $\bigcup_{k=1}^{N} \mathcal{E} ( S_k,  N\Z {+} \rho(k) )$ with $\rho \in \boldsymbol{\mathrm{S}}_N$ cannot be a Riesz basis for $L^2(\bigcup_{k=1}^{N} S_k)$. See Section~\ref{sec:discussion-and-outlook} for some necessary conditions for $\bigcup_{k=1}^K \mathcal{E}(S_k, \Lambda_k)$ to be a Riesz basis for $L^2( \bigcup_{k=1}^K S_k )$.

\begin{corollary}\label{cor:rational-endpoints}
Let $S_1, \ldots, S_{K} \subset [0,1)$ be finite unions of intervals with rational endpoints.
There exist pairwise disjoint sets $\Lambda_1, \ldots, \Lambda_{K} \subset \Z$ such that
\[
		\bigcup_{k=1}^{K} \mathcal{E}(S_k, \Lambda_k)=\bigcup_{k=1}^{K}\big\{ e^{2\pi i \lambda\, (\cdot)}\chi_{S_k} : \lambda\in \Lambda_k\big\}
\]
is a Riesz basis for $L^2( \bigcup_{k=1}^{K} S_k )$.
In particular, if $\bigcup_{k=1}^{K} S_k = [0,1)$, then $\Lambda_1, \ldots, \Lambda_{K}$ form a partition of $\Z$.
\end{corollary}

Theorem~\ref{thm:rational-endpoints} relies on the following two auxiliary results.

\begin{lemma}\label{lem:permuation-exists} 
	For any $A\in\mathbb C^{K\times K}$ and a binary matrix $M\in\{0,1\}^{K\times K}$, there exists a permutation $\rho \in \boldsymbol{\mathrm{S}}_K$ with
\[
\big| \det ((P_\rho A)\odot M) \big|  \geq  \frac{R \cdot \det (A)}{K!} ,
\]
where $P_\rho$ is the permutation matrix associated with $\rho$ (more precisely, $P_\rho = [ p_{k,\ell} ]_{k,\ell=1}^{K}$ with $p_{k,\ell} = 1$ if $\ell = \rho(k)$ and $p_{k,\ell} = 0$ otherwise),
the symbol ${\odot}$ denotes the entrywise multiplication (Hadamard product), and $R$ denotes the number of distinct generalized diagonals\footnote{A generalized diagonal of a $K {\times} K$ matrix is a selection of $K$ cells from the $K {\times} K$ cells, such that exactly one cell is selected from each row and each column. For instance, the usual diagonal is a generalized diagonal.} of $M$ that contain no zeros.

Consequently, if $A\in\mathbb C^{K\times K}$ is nonsingular and if $M\in\{0,1\}^{K\times K}$ has at least one generalized diagonal of ones, then there exists a permutation matrix $P$ with $\det ((PA) \odot M)\neq 0$.
\end{lemma}

In the proof of Theorem~\ref{thm:rational-endpoints}, we will apply Lemma~\ref{lem:permuation-exists} to square submatrices $A$ of the Fourier matrix $W_N = [ e^{- 2\pi i k \ell /N } ]_{k , \ell \in \Z_N}$.

\begin{remark}
\rm
Lemma~\ref{lem:permuation-exists} implies that for any invertible matrix $A\in\mathbb C^{N\times N}$ and a mask $M\in\{0,1\}^{N\times N}$ with at least one generalized diagonal of ones, there exists a row permutation of $M$, say $\widetilde{M}$, such that $A \odot \widetilde{M}$ is invertible.
\end{remark}

The following proposition provides a simple characterization for $\bigcup_{k=1}^{K}\mathcal{E}(S_k, \Lambda_k)$ to be a Riesz basis for $L^2( \bigcup_{k=1}^{K} S_k )$
when $\Lambda_k$ are $N$-periodic sets in $\R$ and $S_k$ are some unions of $[\frac{\ell}{N},\frac{\ell+1}{N})$, $\ell \in \Z_N := \{0, \ldots, N{-}1 \}$.
For notational convenience, we will use indices $k=0,\ldots,K{-}1$, instead of $k=1,\ldots,K$. 

\begin{proposition}\label{prop:characterization-Fourier-submatrices-combination} 
Fix any $N \in \N$.
Let $c_0, \ldots , c_{K-1} \in [0,N)$ be distinct real numbers and let $\mathcal{L}_k \subset \Z_N$, $k=0,\ldots,K{-}1$.
Define $\mathcal{L} := \bigcup_{k=0}^{K-1} \mathcal{L}_k$, $L:= |\mathcal{L}|$, and $W_{(c_0,\mathcal{L}_0),\ldots,(c_{K-1},\mathcal{L}_{K-1})} := [ w_{k,\ell} ]_{k=0,\ldots,K-1, \; \ell \in \mathcal{L}}  \in \C^{K \times L}$ with 
\[
w_{k,\ell}
=
\begin{cases}
e^{- 2\pi i c_k \ell / N} & \text{if} \;\; \ell \in \mathcal{L}_k , \\
0 & \text{if} \;\; \ell \in \mathcal{L} \backslash \mathcal{L}_k .
\end{cases}
\]
Also, define $S_k := \bigcup_{\ell \in \mathcal{L}_k} [\frac{\ell}{N},\frac{\ell+1}{N})$ for $k=0,\ldots,K{-}1$, and $S := \bigcup_{k=0}^{K-1} S_k = \bigcup_{\ell \in \mathcal{L}} \, [\frac{\ell}{N},\frac{\ell+1}{N})$. 
Then
\begin{equation}\label{eqn:exp-system-chi-combination}
\bigcup_{k=0}^{K-1} \mathcal{E} ( S_k,  N\Z {+} c_k )
=\bigcup_{k=0}^{K-1}\big\{ e^{2\pi i \lambda\, (\cdot)}\chi_{S_k} : \lambda\in N\Z {+} c_k \big\} 
\end{equation}
is
\begin{itemize}
\item
a frame for $L^2(S)$ if and only if the mapping $x \mapsto W_{(c_0,\mathcal{L}_0),\ldots,(c_{K-1},\mathcal{L}_{K-1})} x$ is injective, i.e., the matrix $W_{(c_0,\mathcal{L}_0),\ldots,(c_{K-1},\mathcal{L}_{K-1})}$ has full rank and $K \geq L$. 

\item
a Riesz sequence in $L^2(S)$ if and only if the mapping $x \mapsto W_{(c_0,\mathcal{L}_0),\ldots,(c_{K-1},\mathcal{L}_{K-1})} x$ is surjective, i.e., the matrix $W_{(c_0,\mathcal{L}_0),\ldots,(c_{K-1},\mathcal{L}_{K-1})}$ has full rank and $K \leq L$. 

\item
a Riesz basis for $L^2(S)$ if and only if the mapping $x \mapsto W_{(c_0,\mathcal{L}_0),\ldots,(c_{K-1},\mathcal{L}_{K-1})} x$ is bijective, i.e., the matrix $W_{(c_0,\mathcal{L}_0),\ldots,(c_{K-1},\mathcal{L}_{K-1})}$ is invertible (and $K=L$). 
\end{itemize}
In any of the above cases, the optimal lower frame/Riesz bound is given by $\frac{1}{N} \, \sigma_{\min}^2 (W_{(c_0,\mathcal{L}_0),\ldots,(c_{K-1},\mathcal{L}_{K-1})})$.
Furthermore, if $K = L = N$ and if the matrix $W_{(c_0,\mathcal{L}_0),\ldots,(c_{K-1},\mathcal{L}_{K-1})}$ is invertible, then the dual Riesz basis of \eqref{eqn:exp-system-chi-combination} in $L^2[0,1)$ is given by
\begin{equation}\label{eqn:exp-system-chi-combination-dualRB}
\begin{split}
&\bigcup_{k=0}^{N-1} \Big( N  \sum_{j=0}^{N-1} e^{- 2\pi i c_k j / N } \, z_{j,k} \,  \chi_{[\frac{j}{N}, \frac{j+1}{N})} \Big) \cdot \mathcal{E} ( [0,1) ,  N\Z {+} c_k ) \\
&=\bigcup_{k=0}^{N-1} \Big\{ N  \sum_{j=0}^{N-1} e^{- 2\pi i c_k j / N } \, z_{j,k} \,  \chi_{[\frac{j}{N}, \frac{j+1}{N})} (\cdot) \, e^{2\pi i \lambda\, (\cdot)}  : \lambda\in N\Z {+} c_k \Big\} ,
\end{split}
\end{equation}
where $\big[ z_{j,k} \big]_{j,k \in \Z_N} = \big( W_{(c_0,\mathcal{L}_0),\ldots,(c_{N-1},\mathcal{L}_{N-1})} \big)^{-1}$.
\end{proposition}

Unfortunately, if some endpoints of $S_k$ are irrational, there is no convenient characterization for $\bigcup_{k=1}^{K}\mathcal{E}(S_k, \Lambda_k)$ to be a Riesz basis for $L^2( \bigcup_{k=1}^{K} S_k )$.
Nevertheless, we are able to prove the following as our second main result.

\begin{theorem}\label{thm:K-intervals-RB} 
Let $I_1, I_2, I_3$ be intervals partitioning $[0,1)$  
and let $\Lambda_1, \Lambda_2, \Lambda_3$ be a partition of $\Z$ such that for each $k=1, 2, 3$, the system $\mathcal{E} (I_k, \Lambda_k)$ is a Riesz basis for $L^2(I_k)$. 
Here, the intervals are allowed to be empty sets, and we set $\Lambda_k = \emptyset$ in the case that $I_k = \emptyset$.
Then
\[
\bigcup_{k=1}^{3} \mathcal{E}(S_k, \Lambda_k)
=\bigcup_{k=1}^{3} \big\{ e^{2\pi i \lambda (\cdot)}\chi_{S_k} : \lambda\in \Lambda_k\big\}
\]
is a Riesz basis for $L^2( \bigcup_{k=1}^{3} S_k ) = L^2[0,1)$ whenever $S_k=\bigcup_{n\in \mathcal{L}_k}I_n$ with $k \in \mathcal{L}_k \subset\{1, 2, 3 \}$ for all $k=1, 2, 3$. 
\end{theorem}

Note that the condition $k \in \mathcal{L}_k$ is equivalent to having $I_k \subset S_k$.
This condition is actually not necessary but is a convenient assumption. We will discuss this in more detail in Section~\ref{sec:discussion-and-outlook}.

It is easily seen that Theorem~\ref{thm:K-intervals-RB} does not hold for more than three intervals (see Example \ref{ex:K4} below).

We would like to highlight that the statements of Theorems~\ref{thm:rational-endpoints} and \ref{thm:K-intervals-RB} have completely different quantifiers.
Given sets $S_1, \ldots, S_K \subset [0,1)$ with certain property,
Theorem~\ref{thm:rational-endpoints} finds some pairwise disjoint sets $\Lambda_1, \ldots, \Lambda_K \subset \Z$ such that $\bigcup_{k=1}^K \mathcal{E}(S_k, \Lambda_k)$ is a Riesz basis for $L^2( \bigcup_{k=1}^K S_k )$.
In contrast, Theorem~\ref{thm:K-intervals-RB} assumes the sets $I_1, \ldots, I_K \subset [0,1)$ and $\Lambda_1, \ldots, \Lambda_K \subset \Z$ to satisfy that $\mathcal{E} (I_k, \Lambda_k)$ is a Riesz basis for $L^2(I_k)$, $k=1, \ldots, K$, and then shows that
$\bigcup_{k=1}^K \mathcal{E}(S_k, \Lambda_k)$ is a Riesz basis for $L^2( \bigcup_{k=1}^K S_k ) = L^2[0,1)$ if each $S_k$ is a union of subcollection of $I_1, \ldots, I_K$, including $I_k$.

\subsection{Related work}
\label{subsec:RelatedWork}

\

Kozma and Nitzan \cite{KN15} proved that for any finite union $S$ of disjoint intervals in $[0,1)$, there exists a set $\Lambda \subset \Z$ such that $\mathcal{E} (S, \Lambda)$ is a Riesz basis for $L^2(S)$.
By adapting the proof technique of \cite{KN15}, Caragea and Lee \cite{CL22} showed that
if $I_k = [a_k,b_k)$, $k=1,\ldots,K$, are disjoint intervals in $[0,1)$ with the numbers $1, a_1, \ldots, a_K, b_1, \ldots, b_K$ being linearly independent over $\Q$,
then there exist pairwise disjoint sets $\Lambda_k \subset \Z$, $k=1, \ldots, K$,
such that for every $\mathcal{K} \subset \{ 1, \ldots , K \}$, the system $\mathcal{E} (\cup_{k \in \mathcal{K}} \, I_k, \cup_{k \in \mathcal{K}} \, \Lambda_k )$ is a Riesz basis for $L^2 ( \cup_{k \in \mathcal{K}} \, I_k)$.
This extends the result of Kozma and Nitzan by requiring the exponential Riesz basis to possess a hierarchical Riesz bases property so that every subcollections of $I_1, \ldots, I_K$ admit Riesz bases with the corresponding frequency sets, when the endpoints of $I_1, \ldots, I_K$ are rationally independent.
Conjecture \ref{conj:rational-endpoints-special-everyK} below states a similar hierarchical result for the case where all the intervals $I_1, \ldots, I_K$ have rational endpoints.

On the other hand, using completely different proof techniques,
Pfander, Revay and Walnut \cite{PRW21} showed that
for any partition of $[0,1)$ into intervals $I_1, \ldots, I_K$, there exists a partition of $\Z$ into sets $\Lambda_1, \ldots, \Lambda_K$, such that
for each $k$, the system $\mathcal{E} (I_k, \Lambda_k)$ is a Riesz basis for $L^2(I_k)$, and moreover for any $\mathcal{K}\subset \{1,\ldots,K\}$, the system
$\bigcup_{k\in \mathcal{K}}\mathcal{E} (I, \Lambda_k)$
is a Riesz basis for $L^2(I)$ if $I \subset \R$ is an interval of length $\sum_{k\in \mathcal{K}} | I_k |$.

\section{Discussion \& Outlook}
\label{sec:discussion-and-outlook}

We first give an example which demonstrates Theorem~\ref{thm:rational-endpoints}.

\begin{example}\label{ex:first-result}
Let $S_1 = S_3 = [0,\frac{1}{4}) \cup [\frac{2}{4},\frac{3}{4})$ and $S_2 = S_4 = [0,1)$, which corresponds to $K = N = 4$ in Theorem~\ref{thm:rational-endpoints}. Using Proposition~\ref{prop:characterization-Fourier-submatrices-combination}, we see that for $\rho \in \boldsymbol{\mathrm{S}}_4 = \boldsymbol{\mathrm{S}} ( \{ 1,2,3,4 \} )$, the system $\bigcup_{k=1}^{4} \mathcal{E} ( S_k,  4\Z {+} \rho(k) )$ is a Riesz basis for $L^2(\bigcup_{k=1}^{4} S_k) = L^2[0,1)$ if and only if the matrix
\[
A(\rho)
:=
\begin{bmatrix}
e^{- 2\pi i \rho(1) \cdot 0 / 4} & 0 & e^{- 2\pi i \rho(1) \cdot 2 / 4} & 0 \\
e^{- 2\pi i \rho(2) \cdot 0 / 4} & e^{- 2\pi i \rho(2) / 4}  & e^{- 2\pi i \rho(2) \cdot 2 / 4}   & e^{- 2\pi i \rho(2) \cdot 3 / 4}   \\
e^{- 2\pi i \rho(3) \cdot 0 / 4} & 0  & e^{- 2\pi i \rho(3) \cdot 2 / 4}   & 0   \\
e^{- 2\pi i \rho(4) \cdot 0 / 4} & e^{- 2\pi i \rho(4) / 4}  & e^{- 2\pi i \rho(4) \cdot 2 / 4}   & e^{- 2\pi i \rho(4) \cdot 3 / 4}
\end{bmatrix}
\]
is invertible. It is easily seen that 
\[
A(\mathrm{id})
=
\begin{bmatrix}
1 & 0 & -1 & 0 \\
1 & -1 & 1 & -1 \\
1 & 0 & -1 & 0 \\
1 & 1 & 1 & 1
\end{bmatrix}
\]
is singular, but
\[
A(\rho)
=
\begin{bmatrix}
1 & 0 & -1 & 0 \\
1 & i & -1 & -i \\
1 & 0 & 1 & 0 \\
1 & 1 & 1 & 1
\end{bmatrix}
\quad \text{for} \;\;
\rho =
\begin{pmatrix}
1 & 2 & 3 & 4\\
1 & 3 & 2 & 4
\end{pmatrix}
\]
is invertible.
\end{example}

Landau's density theorem \cite{La67} states that
if $\mathcal{E} (S,\Lambda)$ is a Riesz basis for $L^2(S)$, then it necessarily holds that $D(\Lambda) = | S |$,
where $D(\cdot)$ denotes the Beurling density\footnote{A discrete set $\Lambda$ has {\it uniform Beurling density} if $D^-(\Lambda) = \liminf_{r \rightarrow \infty} \inf_{x \in \R} |\Lambda \cap [x,x+r) |/r$ and $D^+(\Lambda) = \limsup_{r \rightarrow \infty} \sup_{x \in \R} |\Lambda \cap [x,x+r) |/r$ are equal, and the {\it Beurling density} of $\Lambda$ is defined as $D(\Lambda) := D^-(\Lambda) = D^+(\Lambda)$.} and $| \cdot |$ denotes the Lebesgue measure on $\R$.
Similarly, Proposition~\ref{prop:characterization-Fourier-submatrices-combination} implies that if $\bigcup_{k=1}^{K} \mathcal{E} ( S_k,  N\Z {+} c_k )$ is a Riesz basis for $L^2(\bigcup_{k=1}^{K} S_k) = L^2(S)$ where each $S_k$ is a union of subcollection of $[\frac{\ell}{N},\frac{\ell+1}{N})$, $\ell \in \Z_N$, then $D(\bigcup_{k=1}^{K} (N\Z {+} c_k )) = \frac{K}{N} = | S |$.
This however is a special setting where $\Lambda_k = N\Z {+} c_k$ are cosets of $N\Z$ and $S_k$ are unions of $[\frac{\ell}{N},\frac{\ell+1}{N})$, $\ell \in \Z_N$.

Let us discuss some necessary conditions for $\bigcup_{k=1}^K \mathcal{E}(S_k, \Lambda_k)$ to be a Riesz basis for $L^2( \bigcup_{k=1}^K S_k )$, in a more general setting.
To this end, we will use the setting of Theorem \ref{thm:K-intervals-RB} but with $K$ intervals $I_1, \ldots, I_K$ where $K \in \N$ is arbitrary.
Let $I_1, \ldots, I_K$ be intervals partitioning $[0,1)$ and let $\Lambda_1, \ldots, \Lambda_K$ be a partition of $\Z$ such that
for each $k$ the system $\mathcal{E} (I_k, \Lambda_k)$ is a Riesz basis for $L^2(I_k)$.
If $S_1, \ldots, S_K$ are unions of $I_1, \ldots, I_K$ such that $\bigcup_{k=1}^K \mathcal{E}(S_k, \Lambda_k)$ is a Riesz basis for $L^2( \bigcup_{k=1}^K S_k )$, then it holds necessarily that
\begin{itemize}

\item[]
(NC1)\; 
$| S_k | \geq D(\Lambda_k) = | I_k |$ for all $k=1,\ldots,K$ \; (due to {\it $\omega$-independence}),

\item[]
(NC2)\;
$D\big( \bigcup_{\ell \,:\, S_\ell \supset I_k} \Lambda_\ell \big) \geq |I_k|$ for all $k=1,\ldots,K$ \; (due to {\it completeness}).

\end{itemize}
We give an example where (NC1) and (NC2) are violated, so that $\bigcup_{k=1}^K \mathcal{E}(S_k, \Lambda_k)$ is not a Riesz basis for $L^2( \bigcup_{k=1}^K S_k )$.

\begin{example}\label{ex:Sk-union-insufficient}
Let $N \geq 5$ be an integer.
Let $I_1 = [0, \frac{N-2}{N})$, $I_2 = [\frac{N-2}{N},\frac{N-1}{N})$, $I_3 = [\frac{N-1}{N},1)$, and $\Lambda_1 = \bigcup_{k=0}^{N-3} (N\Z{+}k)$, $\Lambda_2 = N\Z{+}N{-}2$, $\Lambda_3 = N\Z{+}N{-}1$.
Let $S_k = [0,1) \backslash I_k = \bigcup_{\ell \neq k} I_{\ell}$ for $k=1,2,3$.
Then
$\mathcal{E} (S_1, \Lambda_1) \cup \mathcal{E} (S_2, \Lambda_2) \cup \mathcal{E} (S_3, \Lambda_3)$
is neither $\omega$-independent nor complete in $L^2[0,1)$.
Indeed, the subsystem
$\mathcal{E} (S_1, \Lambda_1)
= \mathcal{E} ([\frac{N-2}{N},1) , \cup_{k=0}^{N-3} N\Z{+}k )$
is not $\omega$-independent since $| [\frac{N-2}{N},1) | = \frac{2}{N} < D ( \bigcup_{k=0}^{N-3} (N\Z{+}k) ) = \frac{N-2}{N}$.
Note that for $x \in I_1$,
\[
\begin{split}
\sum_{k=1}^3 \sum_{\lambda\in \Lambda_k} c_\lambda \, e^{2\pi i \lambda x} \chi_{S_k} (x)
&= \sum_{\lambda\in \Lambda_2} c_\lambda \, e^{2\pi i \lambda x}  + \sum_{\lambda\in \Lambda_3} c_\lambda \, e^{2\pi i \lambda x} .
\end{split}
\]
Since $D(\Lambda_2 \cup \Lambda_3) = \frac{2}{N}$, this series expansion cannot express all functions in $L^2(I_1) = L^2[0, \frac{N-2}{N})$, hence, $\mathcal{E} (S_1, \Lambda_1) \cup \mathcal{E} (S_2, \Lambda_2) \cup \mathcal{E} (S_3, \Lambda_3)$ is not complete in $L^2[0,1)$.
\end{example}

We also provide an example showing that (NC1) and (NC2) are not sufficient.

\begin{example}
\label{ex:K4}
Let $I_k=[\frac{k-1}{4},\frac{k}{4})$ and $\Lambda_k = 4\Z{+}(k{-}1)$ for $k=1,2,3,4$.
The system $\mathcal{E} (I_1 \cup I_3, 4\Z \cup (4\Z{+}2) )$ is not a Riesz basis for $L^2(I_1 \cup I_3)$ and similarly, $\mathcal{E} (I_2 \cup I_4, (4\Z{+}1) \cup (4\Z{+}3) )$ is not a Riesz basis for $L^2(I_2 \cup I_4)$.
Since $L^2[0,1) = L^2(I_1 \cup I_3) \oplus L^2(I_2 \cup I_4)$,
it follows that
$\mathcal{E} (I_1 \cup I_3, 4\Z \cup (4\Z{+}2) ) \cup \mathcal{E} (I_2 \cup I_4, (4\Z{+}1) \cup (4\Z{+}3) )$ is not a Riesz basis for $L^2[0,1)$.
However, the sets $S_1 = S_3 = I_1 \cup I_3$ and $S_2 = S_4 = I_2 \cup I_4$ satisfy (NC1) and (NC2).
\end{example}

Note that both (NC1) and (NC2) are fulfilled if $I_k \subset S_k$ for all $k=1,\ldots,K$.
The condition $I_k \subset S_k$ is not necessary; for instance,
$\mathcal{E} (I_2, \Lambda_1) \cup \mathcal{E} (I_3, \Lambda_2) \cup \mathcal{E} (I_1, \Lambda_3)$ is a Riesz basis for $L^2[0,1)$
when $I_k=[\frac{k-1}{3},\frac{k}{3})$ and $\Lambda_k = 3\Z{+}(k{-}1)$ for $k=1,2,3$.
However, it is a convenient assumption for our applications.
This justifies the assumption $k \in \mathcal{L}_k$ ($\Leftrightarrow$ $I_k \subset S_k$) in Theorem~\ref{thm:K-intervals-RB}.

We now state some conjectures motivated by our results.
For notational convenience, we will use indices $k$ starting from $0$, which allows us to relate with the Fourier matrix $W_N = [ e^{- 2\pi i k \ell /N } ]_{k , \ell \in \Z_N}$. 

\begin{Conjecture}
\label{conj:rational-endpoints-special}
Given a prime $N \in \N$, there exists a permutation $\rho \in \boldsymbol{\mathrm{S}} (\Z_N)$ such that
$\bigcup_{k=0}^{N-1}\mathcal{E} (S_k, N\mathbb Z {+} \rho(k) )$
is a Riesz basis for $L^2( \bigcup_{k=0}^{N-1} S_k ) = L^2[0,1)$
whenever $S_k=\bigcup_{\ell \in \mathcal{L}_k} [\frac{\ell}{N}, \frac{\ell+1}{N})$ and $k \in \mathcal{L}_k \subset \{0,\ldots,N{-}1\}$ for all $k=0,\ldots,N{-}1$.

Equivalently, given a prime $N \in \N$, there exists a permutation $\rho \in \boldsymbol{\mathrm{S}} (\Z_N)$ such that for any binary matrix $M \in\{0,1\}^{N\times N}$ with nonvanishing main diagonal,  
the matrix
$(P_\rho W_N) \odot M = [ e^{- 2\pi i \rho(k) \ell /N } ]_{k, \ell \in \Z_N}  \odot M$
is invertible.
\end{Conjecture}

Conjecture \ref{conj:rational-endpoints-special} generalizes Theorem~\ref{thm:rational-endpoints} by allowing $S_k$ to be {\it arbitrary} unions of $[\frac{\ell}{N}, \frac{\ell+1}{N})$, $\ell \in \Z_N$, including $[\frac{k}{N}, \frac{k+1}{N})$ while the sets $\Lambda_k = N\mathbb Z {+} \rho(k)$ are fixed upon the choice of $\rho \in \boldsymbol{\mathrm{S}} (\Z_N)$, but for $N$ prime.

Numerical experiments show that for $N=4$, Conjecture \ref{conj:rational-endpoints-special} holds with the permutation $\rho \in \boldsymbol{\mathrm{S}} (\Z_4)$ given by $\rho(0) = 0$, $\rho(1) = 2$, $\rho(2) = 1$, $\rho(3) = 3$. 
For $N=5$, the conjecture holds with $\rho \in \boldsymbol{\mathrm{S}} (\Z_5)$ given by $\rho(0) = 0$, $\rho(1) = 1$, $\rho(2) = 2$, $\rho(3) = 4$, $\rho(4) = 3$. 
However, it is confirmed numerically that for $N = 6$, there is no permutation $\rho \in \boldsymbol{\mathrm{S}} (\Z_N)$ satisfying the statement of Conjecture \ref{conj:rational-endpoints-special}.
This is the reason for limiting $N \in \N$ to prime numbers in Conjecture \ref{conj:rational-endpoints-special}.

On the other hand, numerical experiments suggest that for any prime $N \neq 5$, Conjecture \ref{conj:rational-endpoints-special} holds with the identity permutation $\rho = \mathrm{id}$,
that is,
the Fourier matrix $W_N = [ e^{- 2\pi i k \ell /N } ]_{k , \ell \in \Z_N}$ with some off-diagonal entries set to $0$, is always invertible.
For $N = 5$, the identity permutation $\rho = \mathrm{id}$ does not work because, for instance, the matrices  
\[
\begin{bmatrix}
1 & 0 & 0 & 0 & 0 \\
0 & \omega & 0 & \omega^3 & 0 \\
0 & \omega^2 & \omega^4 & 0 & 0 \\
0 & 0 & 0 & \omega^4 & \omega^2 \\
0 & 0 & \omega^3 & 0 & \omega
\end{bmatrix}
\quad
\text{and}
\quad
\begin{bmatrix}
1 & 0 & 0 & 0 & 0 \\
0 & \omega & \omega^2 & 0 & 0 \\
0 & 0 & \omega^4 & 0 & \omega^3 \\
0 & \omega^3 & 0 & \omega^4 & 0 \\
0 & 0 & 0 & \omega^2 & \omega
\end{bmatrix}
\quad \text{with} \;\; \omega = e^{-2 \pi i / 5}
\]
are singular.  
Note that for both matrices, the lower right $4 {\times} 4$ submatrix is singular due to a symmetry.

By restricting $M \in\{0,1\}^{N\times N}$ to the class of block diagonal matrices with all diagonal blocks being full of ones,
we obtain a statement weaker than Conjecture \ref{conj:rational-endpoints-special}.
Note, for instance, that if $N=5$ and  
\[
M = 
\begin{bmatrix}
1 & 1 & 0 & 0 & 0 \\
1 & 1 & 0 & 0 & 0 \\
0 & 0 & 1 & 1 & 0 \\
0 & 0 & 1 & 1 & 0 \\
0 & 0 & 0 & 0 & 1
\end{bmatrix} , 
\]
then the matrix $[ e^{- 2\pi i \rho(k) \ell /N } ]_{k, \ell \in \Z_N}  \odot M$
is invertible if and only if $[ e^{- 2\pi i \rho(k) \ell /N } ]_{k, \ell \in \mathcal{K}}$ is invertible for $\mathcal{K} = \{ 0, 1 \}, \{ 2,3 \}, \{ 4 \}$.
We however allow $N \in \N$ to be composite numbers in contrast with Conjecture \ref{conj:rational-endpoints-special}.

\begin{Conjecture}
\label{conj:rational-endpoints-special-everyK}
Given any $N \in \N$, there exists a permutation $\rho \in \boldsymbol{\mathrm{S}} (\Z_N)$ such that for every $\mathcal{K} \subset \{0,\ldots,N{-}1\}$, the system
$\mathcal{E} (S_{\mathcal{K}}, \bigcup_{k \in \mathcal{K}} N\mathbb Z {+} \rho(k))$
is a Riesz basis for $L^2(S_{\mathcal{K}})$ with $S_{\mathcal{K}} = \bigcup_{k \in \mathcal{K}} [\frac{k}{N},\frac{k+1}{N})$.

Equivalently, given any $N \in \N$, there exists a permutation $\rho \in \boldsymbol{\mathrm{S}} (\Z_N)$ such that $[ e^{- 2\pi i \rho(k) \ell /N } ]_{k, \ell \in \mathcal{K}}$ is invertible for every $\mathcal{K} \subset \{0,\ldots,N{-}1\}$, that is, all principal submatrices of $P_\rho W_N = [ e^{- 2\pi i \rho(k) \ell /N } ]_{k, \ell \in \Z_N}$ are invertible.
\end{Conjecture}

Note that Conjecture \ref{conj:rational-endpoints-special-everyK} holds for all primes $N \in \N$,
due to Chebotar\"{e}v's theorem on roots of unity which asserts that every minor of the Fourier matrix in prime dimension is nonzero, see e.g., \cite{SL96}.
Numerical evidence shows that Conjecture \ref{conj:rational-endpoints-special-everyK} holds also for composite numbers $N \in \N$ at least up to $18$.
For instance, Conjecture \ref{conj:rational-endpoints-special-everyK} for $N = 4$ holds with $\rho \in \boldsymbol{\mathrm{S}} (\Z_4)$ given by $\rho(0) = 0$, $\rho(1) = 2$, $\rho(2) = 1$, $\rho(3) = 3$. Note that the identity permutation $\rho = \mathrm{id} \in \boldsymbol{\mathrm{S}} (\Z_4)$ does not work because the matrices 
\[
\big[ e^{- 2\pi i k \ell /4} \big]_{k, \ell \in \{ 0 , 2 \} }
=
\begin{bmatrix}
1 & 1 \\
1 & 1
\end{bmatrix}
\quad
\text{and}
\quad
\big[ e^{- 2\pi i k \ell /4} \big]_{k, \ell \in \{ 1 , 3 \} }
=
\begin{bmatrix}
-i & i \\
i & -i
\end{bmatrix}
\]
are singular.

Resolving Conjecture \ref{conj:rational-endpoints-special-everyK} would establish a hierarchical result similar to \cite{CL22} using integer frequencies but for the case where $I_1, \ldots, I_K \subset [0,1)$ are disjoint intervals with rational endpoints.  
We mention that in this case, it is relatively easy to derive a hierarchical result with {\it non-integer} frequencies.
Indeed, given any $N \in \N$ and any prime $P > N$, it holds that for every $\mathcal{K} \subset \{0,\ldots,N{-}1\}$ the system
$\mathcal{E} (S_{\mathcal{K}}, \bigcup_{k \in \mathcal{K}} N\mathbb Z {+} \frac{kN}{P})$
is a Riesz basis for $L^2(S_{\mathcal{K}})$ with $S_{\mathcal{K}} = \bigcup_{k \in \mathcal{K}} [\frac{k}{N},\frac{k+1}{N})$ (see \cite[Example 7]{Le22} for more details).

\section{Applications in Sampling Theory}
\label{sec:applications}

We now discuss some applications of our results in sampling theory.

Let $PW_{[-\frac{1}{2},\frac{1}{2})} = \{ f \in L^2(\R) : \mathrm{supp} \, \widehat{f} \subset [-\frac{1}{2},\frac{1}{2}) \}$ be the Paley-Wiener space of signals bandlimited to $[-\frac{1}{2},\frac{1}{2})$.
Here, $\widehat{f} (\omega) = \mathcal{F} [f] (\omega) = \int_{-\infty}^{\infty} f(t) \, e^{- 2 \pi i t \omega} \, dt$ denotes the Fourier transform of $f \in L^2(\R)$.
The celebrated Whittaker-Shannon-Kotel'nikov sampling theorem \cite{Hi96} states that
\[
f(t) = \sum_{n \in \Z} f(n) \, \frac{\sin (\pi (t-n))}{\pi (t-n)}
\quad \text{for all} \;\; f \in PW_{[-\frac{1}{2},\frac{1}{2})} ,
\]
which allows one to reconstruct $f \in PW_{[-\frac{1}{2},\frac{1}{2})}$ exactly from its sample values $\{ f (n) \}_{n \in \Z}$.
The proof is based on the fact that $\mathcal{E} ( [-\frac{1}{2},\frac{1}{2}) , \Z)$ is an orthonormal basis for $L^2[-\frac{1}{2},\frac{1}{2})$.

As a generalization, Papoulis \cite{Pa77} introduced a multi-channel sampling theorem for $PW_{[-\frac{1}{2},\frac{1}{2})}$ using generalized samples. 
We will particularly consider channels which are bandpass filters.
For instance, the channel $f \mapsto \phi_\delta * f$ with a low-pass filter $\phi_\delta$ given by $\widehat{\phi}_\delta = \chi_{[-\delta,\delta)}$ for a small $\delta > 0$, produces output signals with low frequency components only.

Let $S_0, \ldots, S_{N-1} \subset [-\frac{1}{2},\frac{1}{2})$ be unions of $I_\ell := -\frac{1}{2} + [\frac{\ell}{N}, \frac{\ell+1}{N})$, $\ell \in \Z_N$, such that $S_k \supset I_k$ for all $k = 0, \ldots, N{-}1$.
Clearly, we have $\bigcup_{k=0}^{N-1} S_k = [-\frac{1}{2},\frac{1}{2})$.  
For $k=0,\ldots,N{-}1$, we define $\phi_{S_k} \in L^2(\R)$ by
$\widehat{\phi}_{S_k} = \chi_{S_k}$.
By Theorem~\ref{thm:rational-endpoints}, there exists a permutation $\rho \in \boldsymbol{\mathrm{S}} (\Z_N)$ such that
$\bigcup_{k=0}^{N-1} \mathcal{E} ( S_k,  N\Z {+} \rho(k) )$
is a Riesz basis for $L^2 [-\frac{1}{2},\frac{1}{2})$. 
Here, we have used that the complex exponentials with integer frequencies are $1$-periodic, which allows us to use the spaces $L^2 [-\frac{1}{2},\frac{1}{2})$ and $L^2 [0,1)$ interchangeably.  
Adapting Equation~\eqref{eqn:dual-expansion} below (see Proposition \ref{prop:characterization-Fourier-submatrices-combination}), we have that for all $f \in PW_{[-\frac{1}{2},\frac{1}{2})}$, 
\[
\widehat{f} (\omega)
= \sum_{k=0}^{N-1} \sum_{m \in \Z}  \big\langle \widehat{f}, e^{2\pi i (- N m + \rho(k) ) \, (\cdot)} \chi_{S_k}  \big\rangle_{L^2 [-\frac{1}{2},\frac{1}{2}) } \, N \, e^{2\pi i (- N m + \rho(k) ) \omega} \, G_k(\omega) ,
\]
where $G_0, \ldots, G_{N-1}$ are some piecewise constant functions supported in $[-\frac{1}{2},\frac{1}{2})$; more precisely, $G_k (\omega) = \sum_{j=0}^{N-1} e^{- 2\pi i \rho(k) j / N } \, z_{j,k} \,  \chi_{-\frac{1}{2} +[\frac{j}{N}, \frac{j+1}{N})} (\omega)$ for some suitable $z_{0,k}, \ldots, z_{N-1,k} \in \C$. Note that
\[
\begin{split}
\big\langle \widehat{f}, e^{2\pi i (- N m + \rho(k) ) \, (\cdot)} \chi_{S_k}  \big\rangle_{L^2 [-\frac{1}{2},\frac{1}{2}) }
&= \int_{-1/2}^{1/2} \widehat{f} (\omega) \chi_{S_k} (\omega) \, e^{2\pi i (N m - \rho(k) ) \omega} \, d\omega \\
&= f * {\phi}_{S_k} (N m {-} \rho(k)) .
\end{split}
\]
Consequently, we obtain the sampling expansion
\[
f(t)
= \sum_{k=0}^{N-1} \sum_{m \in \Z}   f * {\phi}_{S_k} (N m {-} \rho(k))  \;  g_k( t {-} N m {+} \rho(k) ) ,
\quad f \in PW_{[-\frac{1}{2},\frac{1}{2})} ,
\]
where $g_k := N \cdot \mathcal{F}^{-1} [G_k]$ for $k=0,\ldots,N{-}1$.
This provides an explicit reconstruction formula for $f \in PW_{[-\frac{1}{2},\frac{1}{2})}$ using generalized samples
$\{ f * {\phi}_{S_k} (N m {-} \rho(k)) : k=0,\ldots,N{-}1, \; m \in \Z \}$.

\section{Proof of the first main result}
\label{sec:proof-thm:rational-endpoints}

\subsection{Proof of Theorem~\ref{thm:rational-endpoints}}
\label{subsec:proof-thm:rational-endpoints}

\


For each $k=1,\ldots,K$, we may write $S_k = \bigcup_{\ell \in \mathcal{L}_k} \, [\frac{\ell}{N},\frac{\ell+1}{N})$ for some nonempty set $\mathcal{L}_k \subset \Z_N$.
Since $\bigcup_{k=1}^{K} S_k = \bigcup_{k=1}^{K} I_k$, the set $\mathcal{L} := \bigcup_{k=1}^{K} \mathcal{L}_k$ has cardinality $| \mathcal{L} | = K$.
Let $A =  [ e^{- 2\pi i k \ell / N} ]_{1 \leq k \leq K, \; \ell \in \mathcal{L}}$ and $M = [ m_{k,\ell} ]_{1 \leq k \leq K, \; \ell \in \mathcal{L}}\in \{ 0, 1 \}^{K \times K}$ with
\[
m_{k,\ell}
=
\begin{cases}
1 & \text{if} \;\; \ell \in \mathcal{L}_k , \\
0 & \text{if} \;\; \ell \in \mathcal{L} \backslash \mathcal{L}_k .
\end{cases}
\]
Note that $A$ is an {\it invertible} Vandermonde matrix and by definition of $S_1, \ldots, S_K$, the matrix $M$ has at least one generalized diagonal of ones.
Applying Lemma~\ref{lem:permuation-exists}, we obtain a permutation $\rho \in \boldsymbol{\mathrm{S}}_K$ such that $([ e^{- 2\pi i \rho(k) \ell / N} ]_{1 \leq k \leq K, \; \ell \in \mathcal{L}}) \odot M$ is invertible.
By Proposition~\ref{prop:characterization-Fourier-submatrices-combination}, this means that
$\bigcup_{k=1}^{K}\mathcal{E} ( S_k,  N\Z {+} \rho(k) )$ is a Riesz basis for $L^2(\bigcup_{k=1}^{K} S_k)$.
%
%

\subsection{Proof of Corollary~\ref{cor:rational-endpoints}}
\label{subsec:proof-cor:rational-endpoints}

\

We may assume that $S_1, \ldots, S_{K}$ are unions of $[\frac{\ell}{N}, \frac{\ell+1}{N})$, $\ell = 0, \ldots, N{-}1$, for some $N \in \N$. 
Then $\bigcup_{k=1}^{K} S_k$ is the union of $K' \; ( \leq N)$ distinct intervals of the form $[\frac{\ell}{N}, \frac{\ell+1}{N})$, which are denoted by $I_1, \ldots, I_{K'}$. 
We can choose sets $T_1, \ldots, T_{K'}$ from $S_1, \ldots, S_{K}$ with repetition allowed, such that $I_{k'} \subset T_{k'}$ for all $k' = 1,\ldots,K'$.
It then follows from Theorem~\ref{thm:rational-endpoints} that 
there is a permutation $\rho \in \boldsymbol{\mathrm{S}}_{K'}$ such that
$\bigcup_{k'=1}^{K'} \mathcal{E} ( T_{k'},  N\Z {+} \rho(k') )$ is a Riesz basis for $L^2(\bigcup_{k'=1}^{K'} T_{k'}) = L^2(\bigcup_{k=1}^{K} S_k)$.
For $k=1, \ldots, K$, we define $\Lambda_k$ to be the union of $N\Z {+} \rho(k')$ over all $k'$ such that $T_{k'} = S_k$. Then $\bigcup_{k=1}^{K} \mathcal{E}(S_k, \Lambda_k) = \bigcup_{k'=1}^{K'} \mathcal{E} ( T_{k'},  N\Z {+} \rho(k') )$ is a Riesz basis for $L^2(\bigcup_{k=1}^{K} S_k)$ as desired.

If $\bigcup_{k=1}^{K} S_k = [0,1)$, then $K' = N$ and it follows immediately that $\Lambda_1, \ldots, \Lambda_{K}$ is a partition of $\Z$.

\subsection{Proof of Lemma~\ref{lem:permuation-exists}}
\label{subsec:proof-lem:permuation-exists}

\

Let $A = [ a_{k,\ell} ]_{k,\ell=1}^{K} \in \C^{K \times K}$ and $M = [ m_{k,\ell} ]_{k,\ell=1}^{K} \in \{ 0, 1 \}^{K \times K}$.
Let $G=G_M$ denote the set\footnote{If $M$ contains no zeros, then $G=\boldsymbol{\mathrm{S}}_K$. If $M$ contains a zero row or column, then $G=\emptyset$.} of all permutations $\sigma \in \boldsymbol{\mathrm{S}}_K$ with the property that $m_{k,\sigma(k)}=1$ for $k=1,\ldots,K$.
Note that for any $\rho\in \boldsymbol{\mathrm{S}}_K$, we have $P_{\rho} A = [ a_{\rho(k),\ell} ]_{k,\ell=1}^{K}$.
For $\rho\in \boldsymbol{\mathrm{S}}_K$, we compute
\begin{align*}
	\det ((P_\rho A)\odot M)
	&= \sum_{\sigma\in \boldsymbol{\mathrm{S}}_K} {\rm sgn}(\sigma)\ \prod_{k=1}^{K} a_{\rho(k), \sigma(k)} \, m_{k,\sigma(k)}\\
	&= \sum_{\sigma\in G} {\rm sgn}(\sigma)\ \prod_{k=1}^{K} a_{\rho(k), \sigma(k)} \\
	&= \sum_{\sigma\in G} {\rm sgn}(\sigma)\ \prod_{k=1}^{K} a_{k, (\sigma\circ \rho^{-1}( k))}\\
	&= \sum_{\sigma' \in T_{\rho} G} {\rm sgn}(\sigma' \circ \rho)\ \prod_{k=1}^{K} a_{k, \sigma' ( k)}\\
	&= {\rm sgn}(\rho) \sum_{\sigma\in T_{\rho} G} {\rm sgn}(\sigma)\ \prod_{k=1}^{K} a_{k, \sigma( k)}  , 
\end{align*}
where $\sgn : \boldsymbol{\mathrm{S}}_K \rightarrow \{ \pm 1 \}$ is the sign function of permutations in $\boldsymbol{\mathrm{S}}_K$ (i.e., $\sgn(\sigma) = 1$ if $\sigma$ is an even permutation and $\sgn(\sigma) = -1$ if $\sigma$ is an odd permutation)
and $T_\rho : \boldsymbol{\mathrm{S}}_K \rightarrow \boldsymbol{\mathrm{S}}_K$ is the bijective map given by $T_\rho (\sigma) = \sigma \circ \rho^{-1}$.
Since
\begin{align*}
\sum_{\rho\in \boldsymbol{\mathrm{S}}_K} \Big( \sum_{\sigma\in T_{\rho} G} {\rm sgn}(\sigma)\ \prod_{k=1}^{K} a_{k, \sigma ( k)} \Big) =|G| \cdot \det (A) ,
\end{align*}
there exists at least one $\rho\in \boldsymbol{\mathrm{S}}_K$ with
\begin{align*}
\big| \det ((P_\rho A)\odot M) \big| = \Big|
	 \sum_{\sigma\in T_{\rho} G} {\rm sgn}(\sigma)\ \prod_{k=1}^{K} a_{k, \sigma ( k)} \Big|  \geq  \frac{|G| \cdot \det (A)}{K!} .
\end{align*}
%

\subsection{Proof of Proposition~\ref{prop:characterization-Fourier-submatrices-combination}}
\label{subsec:proof-prop:characterization-Fourier-submatrices-combination}

\

Recall that a sequence $\{ f_n \}_{n\in\Z}$ in a separable Hilbert space $\mathcal{H}$ is called a {\it frame} if there exist constants $0 < A \leq B < \infty$ such that
$A \, \| f \|^2
\leq \sum_{n\in\Z} |\langle f , f_n\rangle|^{2}
\leq B \, \| f \|^2$ for all $f \in \mathcal{H}$.
Since frames are necessarily complete, it follows that $\{ f_n \}_{n\in\Z}$ is a Riesz basis for $\mathcal{H}$ if and only if it is both a frame and a Riesz sequence.


\medskip

\noindent
\textbf{Frame}. \
Fix any $f \in L^2(S)$.
For any $k \in \{ 0, \ldots, K{-}1 \}$ and $\lambda \in N\Z {+} c_{k}$, we have
\[
\begin{split}
\langle f, e^{2\pi i \lambda (\cdot)} \, \chi_{\bigcup_{\ell \in \mathcal{L}_k} [\frac{\ell}{N},\frac{\ell+1}{N})}  \rangle_{L^2(S)}
&= \int_{\bigcup_{\ell \in \mathcal{L}_k} [\frac{\ell}{N},\frac{\ell+1}{N})} f(t) \, e^{- 2\pi i \lambda t} \, dt \\
&= \int_0^{\frac{1}{N}} \sum_{\ell \in \mathcal{L}_k} f ( t {+} \tfrac{\ell}{N} ) \, e^{- 2\pi i \lambda (t + \frac{\ell}{N}) } \, dt \\
&= \int_0^{\frac{1}{N}} F_{k} (t) \, e^{- 2\pi i \lambda t} \, dt
= \langle F_{k} , e^{2\pi i \lambda (\cdot)} \rangle_{L^2[0,\frac{1}{N})}
\end{split}
\]
where
\begin{equation}\label{eqn:def-Fm}
F_{k} (t)
:= \sum_{\ell \in \mathcal{L}_k} f (t {+} \tfrac{\ell}{N} ) \, e^{- 2\pi i c_{k} \ell / N}
\quad \text{for} \;\; t \in [0,\tfrac{1}{N}) .
\end{equation}
Then
\[
\begin{split}
&\sum_{k=0}^{K-1} \sum_{\lambda \in N\Z + c_{k}}
\big| \langle f, e^{2\pi i \lambda (\cdot)} \, \chi_{\bigcup_{\ell \in \mathcal{L}_k} [\frac{\ell}{N},\frac{\ell+1}{N})}  \rangle_{L^2(S)} \big|^2 \\
&=  \sum_{k=0}^{K-1} \sum_{\lambda \in N\Z + c_{k}} \big| \langle F_{k} , e^{2\pi i \lambda (\cdot)} \rangle_{L^2[0,\frac{1}{N})}  \big|^2 
= \tfrac{1}{N} \,  \sum_{k=0}^{K-1}  \| F_{k} \|_{L^2[0,\frac{1}{N})}^2
\end{split}
\]
by using that $\mathcal{E} ( [0,\frac{1}{N}),  N\Z )$ is an orthogonal basis for $L^2[0,\frac{1}{N})$ with tight Riesz bound $\frac{1}{N}$.
Rewriting \eqref{eqn:def-Fm} as
\[
F_{k} (t)
= \sum_{\ell \in \mathcal{L}_k} f (t {+} \tfrac{\ell}{N} ) \cdot e^{- 2\pi i c_{k} \ell / N}
+ \sum_{\ell \in \mathcal{L} \backslash \mathcal{L}_k} f (t {+} \tfrac{\ell}{N} ) \cdot 0
\]
and collecting the equation for $k=0,\ldots,K{-}1$, we get
\begin{equation}\label{eqn:pf-frame-rel}
\big[
F_{k} (t)
\big]_{k=0}^{K-1}
=
W_{(c_0,\mathcal{L}_0),\ldots,(c_{K-1},\mathcal{L}_{K-1})} \,
\big[ f ( t {+} \tfrac{\ell}{N} )
\big]_{\ell \in \mathcal{L}}
\quad \text{for} \;\; t \in [0,\tfrac{1}{N}) .
\end{equation}
(i) If $K < L$, then the kernel of $W_{(c_0,\mathcal{L}_0),\ldots,(c_{K-1},\mathcal{L}_{K-1})} \in \C^{K \times L}$ has dimension at least $L {-} K > 0$.
Fix any nontrivial vector $\bv \in \kernel ( W_{(c_0,\mathcal{L}_0),\ldots,(c_{K-1},\mathcal{L}_{K-1})} ) \subset \C^{L}$ and set $f \in L^2(S)$ by $[ f ( t {+} \frac{\ell}{N} ) ]_{\ell \in \mathcal{L}} = \bv$ for $t \in [0,\frac{1}{N})$, so that $\sum_{k=1}^K \| F_{k} \|_{L^2[0,\frac{1}{N})}^2 = 0$ by \eqref{eqn:pf-frame-rel} while $\| f \|_{L^2(S)}^2 = \frac{1}{N} \cdot \| \bv \|_2^2 \neq 0$.
Hence, the system \eqref{eqn:exp-system-chi-combination} is not a frame for $L^2(S)$. Note that since $K < L$, the mapping $\C^L \rightarrow \C^K$, $x \mapsto W_{(c_0,\mathcal{L}_0),\ldots,(c_{K-1},\mathcal{L}_{K-1})} x$ is not injective. \\
(ii) If $K \geq L$, then
\[
\sum_{k=0}^{K-1} | F_{k} (t) |^2
 \overset{\eqref{eqn:pf-frame-rel}}{\geq} \sigma_{\min}^2 ( W_{(c_0,\mathcal{L}_0),\ldots,(c_{K-1},\mathcal{L}_{K-1})} )  \cdot \sum_{\ell \in \mathcal{L}} | f ( t {+} \tfrac{\ell}{N} ) |^2
\quad \text{for} \;\; t \in [0,\tfrac{1}{N}) ,
\]
and integrating both sides with respect to $t \in [0,\frac{1}{N})$ gives
\[
\sum_{k=0}^{K-1} \| F_{k} \|_{L^2[0,\frac{1}{N})}^2
\;\geq\;
\sigma_{\min}^2 (W_{(c_0,\mathcal{L}_0),\ldots,(c_{K-1},\mathcal{L}_{K-1})})   \cdot \| f \|_{L^2(S)}^2
\]
and therefore
\[
\begin{split}
&\sum_{k=0}^{K-1} \sum_{\lambda \in N\Z + c_{k}}
\big| \langle f, e^{2\pi i \lambda (\cdot)} \, \chi_{\bigcup_{\ell \in \mathcal{L}_k} [\frac{\ell}{N},\frac{\ell+1}{N})}  \rangle_{L^2(S)} \big|^2 \\
&\quad \geq \tfrac{1}{N} \, \sigma_{\min}^2 (W_{(c_0,\mathcal{L}_0),\ldots,(c_{K-1},\mathcal{L}_{K-1})})   \cdot \| f \|_{L^2(S)}^2 .
\end{split}
\]
This shows that if $K \geq L$, then \eqref{eqn:exp-system-chi-combination} is a frame for $L^2(S)$ if and only if the (tall) matrix $W_{(c_0,\mathcal{L}_0),\ldots,(c_{K-1},\mathcal{L}_{K-1})} \in \C^{K \times L}$ has full rank.

\medskip

\noindent
\textbf{Riesz sequence}. \
For any $\boldb = \{ b_\lambda \} \in \ell_2 (\bigcup_{k=0}^{K-1} N\Z {+} c_{k})$, we set
\[
f =
\sum_{k=0}^{K-1} \sum_{\lambda \in N\Z + c_{k}}
b_\lambda \, e^{2\pi i \lambda (\cdot)} \, \chi_{\bigcup_{\ell \in \mathcal{L}_k} [\frac{\ell}{N},\frac{\ell+1}{N})}
\]
which is clearly supported in $S = \bigcup_{\ell \in \mathcal{L}} \, [\frac{\ell}{N},\frac{\ell+1}{N})$.
Then for any $t \in [0,\frac{1}{N})$ and any $\ell \in \Z_N$, we have
\[
f ( t {+} \tfrac{\ell}{N} )
=
\sum_{k=0}^{K-1} \Big( e^{2\pi i c_{k} \ell / N} \, \delta_{\{ \ell \in \mathcal{L}_k \}} \sum_{\lambda \in N\Z + c_{k}} b_\lambda \, e^{2\pi i \lambda t} \Big)
\]
where $\delta_{\{ \ell \in \mathcal{L}_k \}}$ is $1$ if $\ell \in \mathcal{L}_k$ and $0$ otherwise,
so that
\begin{equation}\label{eqn:pf-RS-rel}
\big[
f ( t {+} \tfrac{\ell}{N} )
\big]_{\ell \in \mathcal{L}}
= \big( W_{(c_0,\mathcal{L}_0),\ldots,(c_{K-1},\mathcal{L}_{K-1})} \big)^*
 \big[ \textstyle \sum_{\lambda \in N\Z + c_{k}} b_\lambda \, e^{2\pi i \lambda t}
 \big]_{k=0}^{K-1}  .
\end{equation}
(i) If $K > L$, then the kernel of $( W_{(c_0,\mathcal{L}_0),\ldots,(c_{K-1},\mathcal{L}_{K-1})} )^* \in \C^{L \times K}$ has dimension at least $K {-} L > 0$.
Fix any nontrivial vector $\bv = \{ v_k \}_{k=1}^K \in \kernel ( ( W_{(c_0,\mathcal{L}_0),\ldots,(c_{K-1},\mathcal{L}_{K-1})} )^*  ) \subset \C^K$.
For each $k = 0, \ldots, K{-}1$, the system $\mathcal{E} ( [0,\frac{1}{N}), N\Z {+} c_k )$ is a Riesz basis for $L^2[0,\frac{1}{N})$, so there exists a unique sequence $\{ b_\lambda \} \in \ell_2 (N\Z {+} c_k)$ such that
$\sum_{\lambda \in N\Z + c_k} b_\lambda \, e^{2\pi i \lambda t} = v_k$ for a.e.~$t \in [0,\frac{1}{N})$.
We thus obtain a nontrivial sequence $\boldb = \{ b_\lambda \} \in \ell_2 (\cup_{k=1}^K N\Z {+} c_k)$ satisfying
\[
\big[ \textstyle \sum_{\lambda \in N\Z + c_k} b_\lambda \, e^{2\pi i \lambda t}
\big]_{k=1}^K
=
\bv
\quad \text{for a.e.} \;\; t \in [0,\tfrac{1}{N}) .
\]
The corresponding function $f$ satisfies $\| f \|_{L^2(S)}^2 = \int_0^{1/N} \sum_{\ell \in \mathcal{L}}  | f ( t {+} \frac{\ell}{N} )|^2 \, dt = 0$ by \eqref{eqn:pf-RS-rel} while $\| \boldb \|_2 \neq 0$.
Hence, the system \eqref{eqn:exp-system-chi-combination} is not a Riesz sequence in  $L^2(S)$.
Note that since $K > L$, the mapping $\C^L \rightarrow \C^K$, $x \mapsto W_{(c_0,\mathcal{L}_0),\ldots,(c_{K-1},\mathcal{L}_{K-1})} x$ is not surjective. \\
(ii) If $K \leq L$, then it holds for all $t \in [0,\tfrac{1}{N})$, 
\[
\begin{split}
\sum_{\ell \in \mathcal{L}}  \big| f ( t {+} \tfrac{\ell}{N} ) \big|^2
 \overset{\eqref{eqn:pf-RS-rel}}{\geq}
\sigma_{\min}^2 \big( W_{(c_0,\mathcal{L}_0),\ldots,(c_{K-1},\mathcal{L}_{K-1})} \big)
\cdot \sum_{k=0}^{K-1} \Big| \sum_{\lambda \in N\Z + c_{k}} b_\lambda \, e^{2\pi i \lambda t}  \Big|^2 , 
\end{split}
\]
where we used that $\sigma_{\min}^2 \big(  ( W_{(c_0,\mathcal{L}_0),\ldots,(c_{K-1},\mathcal{L}_{K-1})} )^* \big) = \sigma_{\min}^2 \big( W_{(c_0,\mathcal{L}_0),\ldots,(c_{K-1},\mathcal{L}_{K-1})} \big)$.
Further, integrating both sides with respect to $t \in [0,\frac{1}{N})$ gives
\[
\begin{split}
\| f \|^2_{L^2(S)}
&\geq
\tfrac{1}{N} \, \sigma_{\min}^2 \big( W_{(c_0,\mathcal{L}_0),\ldots,(c_{K-1},\mathcal{L}_{K-1})} \big)   \cdot \| \boldb \|_{\ell_2}^2 ,
\end{split}
\]
where we used that $\mathcal{E} ( [0,\frac{1}{N}),  N\Z )$ is an orthogonal basis for $L^2[0,\frac{1}{N})$ with tight Riesz bound $\frac{1}{N}$.
This shows that if $K \leq L$, then \eqref{eqn:exp-system-chi-combination} is a Riesz sequence in $L^2(S)$ if and only if the (short and fat) matrix $W_{(c_0,\mathcal{L}_0),\ldots,(c_{K-1},\mathcal{L}_{K-1})} \in \C^{K \times L}$ has full rank.
%

\medskip

\noindent
\textbf{Riesz basis}. \
This part follows immediately by combining the {\it frame} and {\it Riesz sequence} parts.

\medskip

\noindent
\textbf{The dual Riesz basis of \eqref{eqn:exp-system-chi-combination} when $K = L = N$}. \
Since both \eqref{eqn:exp-system-chi-combination} and \eqref{eqn:exp-system-chi-combination-dualRB} are Bessel sequences in $L^2[0,1)$, it is enough to show (see e.g., \cite[Lemma 6.3.2]{Ch16}) that for all $f \in L^2[0,1)$,
\begin{equation}\label{eqn:dual-expansion}
\begin{split}
f(t) &= \sum_{k=0}^{N-1} \sum_{m \in \Z} \big\langle f , e^{2\pi i (N m + c_k) \, (\cdot)} \chi_{S_k}  \big\rangle_{L^2[0,1)} \\
&\qquad\qquad\quad \cdot N \, e^{2\pi i (N m + c_k) t}  \sum_{j=0}^{N-1} e^{- 2\pi i c_k j / N } \, z_{j,k} \,  \chi_{[\frac{j}{N}, \frac{j+1}{N})} (t) .
\end{split}
\end{equation}
To see this, note that for any $t \in [\frac{j}{N}, \frac{j+1}{N})$ with $j \in \Z_N$ the right hand side becomes
\[
\begin{split}
&\sum_{k=0}^{N-1} \Big( \sum_{m \in \Z} \Big\langle \sum_{\ell=0}^{N-1} f(\cdot + \tfrac{\ell}{N} ) \, e^{- 2\pi i c_k \, (\cdot + \frac{\ell}{N} )} \chi_{S_k} (\cdot + \tfrac{\ell}{N} )  , e^{2\pi i N m \, (\cdot)} \Big\rangle_{L^2[0,\frac{1}{N})}  \\
&\qquad\qquad\quad \cdot  N \, e^{2\pi i N m t} \Big) z_{j,k} \, e^{2\pi i c_k \, (t - \frac{j}{N} )}  \\
&= \sum_{k=0}^{N-1}  \Big( \sum_{\ell=0}^{N-1} f(t{+} \tfrac{\ell}{N} ) \, e^{- 2\pi i c_k (t + \frac{\ell}{N} )} \chi_{S_k} (t{+} \tfrac{\ell}{N} ) \Big)  z_{j,k} \, e^{2\pi i c_k \, (t - \frac{j}{N} )}  \\
&= \sum_{\ell=0}^{N-1} \Big( \sum_{k=0}^{N-1}  z_{j,k}  \, e^{- 2\pi i c_k (\ell+j)/N } \chi_{S_k} (t{+} \tfrac{\ell}{N} ) \Big)  f(t{+} \tfrac{\ell}{N} )  \\
&= \sum_{\ell=0}^{N-1} \Big( \sum_{k=0}^{N-1}  z_{j,k}  \, w_{k,\ell + j}\Big)  f(t{+} \tfrac{\ell}{N} )  \\
&= \sum_{\ell=0}^{N-1} \delta_{j,\ell + j} \, f(t{+} \tfrac{\ell}{N} )  \\
&= f(t) .
\end{split}
\]
Hence, the dual Riesz basis of \eqref{eqn:exp-system-chi-combination} is given by \eqref{eqn:exp-system-chi-combination-dualRB}.

\section{Proof of Theorem~\ref{thm:K-intervals-RB}}
\label{sec:Proof-Thm-K-intervals-RB}

\

For a sequence $\{ f_n \}_{n\in\Z}$ in a separable Hilbert space $\mathcal{H}$, its {\it synthesis operator} $T : \ell_2(\Z) \rightarrow \mathcal{H}$ is defined by $T( \{ c_n \}_{n\in\Z} ) = \sum_{n\in\Z} c_n \, f_n$.
This operator is well-defined if and only if $\{ f_n \}_{n\in\Z}$ is a Bessel sequence in $\mathcal{H}$ (see e.g., \cite[Theorem~3.2.3, Corollary~3.2.4]{Ch16}).
It is easily seen that $T$ is bijective if and only if $\{ f_n \}_{n\in\Z}$ is a Riesz basis for $\mathcal{H}$.

\begin{theorem}[The Paley-Wiener stability theorem \cite{PW34}; see {\cite[p.\,35, Theorem~13]{Yo01}}]
\label{thm:PaleyWienerStabilityGeneral}
Let $\{ e_n \}_{n\in\Z}$  be an orthonormal basis for a Hilbert space $\mathcal{H}$ and let $\{ f_n \}_{n\in\Z}$  be a sequence in $\mathcal{H}$.
If there is a constant $0 \leq \lambda < 1$ satisfying
\[
\Big\| \sum_{n \in \Z} c_n \, (e_n - f_n) \Big\|_{\mathcal{H}}
\leq \lambda \,  \| \bc \|_{\ell_2}
\quad \text{for all} \;\; \bc = \{ c_n \}_{n \in \Z} \in \ell_2(\Z) ,
\]
then $\{ f_n : n \in \Z \}$ is a Riesz basis for $\mathcal{H}$.
\end{theorem}

We are now ready to prove Theorem~\ref{thm:K-intervals-RB}.
We will denote by $T$ the synthesis operator for $\bigcup_{k=1}^3 \mathcal{E}(S_k, \Lambda_k)$.

\begin{proof}[\textbf{Proof}]
First, if $I_1 = [0,1)$ and $I_2 = I_3 = \emptyset$, then the statement of Theorem~\ref{thm:K-intervals-RB} holds trivially. 

\medskip

Now, we assume that $I_1$ and $I_2$ are {\it nonempty} intervals partitioning $[0,1)$ and $I_3 = \emptyset$. 
There are $2^{1 \cdot 2} = 4$ cases to consider.

\begin{itemize}
\item[Case 1.]
$\mathcal{L}_1 = \{ 1 \}$, $\mathcal{L}_2 = \{ 2 \}$.

Using $L^2[0,1) = L^2(I_1) \oplus L^2(I_2)$, one can easily see that $\mathcal{E} (I_1, \Lambda_1) \cup \mathcal{E} (I_2, \Lambda_2)$ is a Riesz basis for $L^2[0,1)$.

\item[Case 2.]
$\mathcal{L}_1 = \{ 1 \}$, $\mathcal{L}_2 = \{ 1,2 \}$.

\medskip

\noindent
{\it Injectivity of $T$.} \
Let $\bc = \{ c_n \}_{n \in \Z} \in \ell_2(\Z)$ be such that
\[
\sum_{\lambda\in \Lambda_1} c_\lambda \, e^{2\pi i \lambda x} \chi_{I_1} (x) + \sum_{\lambda\in \Lambda_2} c_\lambda \, e^{2\pi i \lambda x} \chi_{I_1 \cup I_2} (x) = 0 \quad
\text{for a.e.~$x \in [0,1)$.  }
\]
Since $\mathcal{E} (I_2, \Lambda_2)$ is a Riesz basis for $L^2(I_2)$, we have $c_\lambda = 0$ for $\lambda \in \Lambda_2$.
In turn, since $\mathcal{E} (I_1, \Lambda_1)$ is a Riesz basis for $L^2(I_1)$, we get $c_\lambda = 0$ for $\lambda \in \Lambda_1$.

\medskip

\noindent
{\it Surjectivity of $T$.} \
Fix any function $f \in L^2[0,1)$.
Since $\mathcal{E} (I_2, \Lambda_2)$ is a Riesz basis for $L^2(I_2)$, there exists a unique sequencs $\{ c_\lambda \}_{\lambda \in \Lambda_2} \in \ell_2(\Lambda_2)$ such that $\sum_{\lambda\in \Lambda_2} c_\lambda \, e^{2\pi i \lambda x} = f(x)$ for a.e.~$x \in I_2$.
Define $f_1 \in L^2(I_1)$ by $f_1(x) := f(x) -  \sum_{\lambda\in \Lambda_2} c_\lambda \, e^{2\pi i \lambda x}$ for $x \in I_1$.
Since $\mathcal{E} (I_1, \Lambda_1)$ is a Riesz basis for $L^2(I_1)$, there exists a unique sequence $\{ c_\lambda \}_{\lambda \in \Lambda_1} \in \ell_2(\Lambda_1)$ such that $\sum_{\lambda\in \Lambda_1} c_\lambda \, e^{2\pi i \lambda x} (x) = f_1(x)$ for a.e.~$x \in I_1$.
Then
\[
\sum_{\lambda\in \Lambda_1} c_\lambda \, e^{2\pi i \lambda x} \chi_{I_1} (x) + \sum_{\lambda\in \Lambda_2} c_\lambda \, e^{2\pi i \lambda x} \chi_{I_1 \cup I_2} (x)  = f(x) \quad
\text{for a.e.~$x \in [0,1)$.  }
\]

Hence, the system $\mathcal{E} (I_1 \cup I_2, \Lambda_1) \cup \mathcal{E}(I_2, \Lambda_2)$ is a Riesz basis for $L^2[0,1)$.

\item[Case 3.]
$\mathcal{L}_1 = \{ 1,2 \}$, $\mathcal{L}_2 = \{ 2 \}$.
Symmetric with Case 2.

\item[Case 4.]
$\mathcal{L}_1 = \{ 1,2 \}$, $\mathcal{L}_2 = \{ 1,2 \}$.
This case is trivial, since $\mathcal{E} ([0,1) , \Z)$ is an orthonormal basis for $L^2[0,1)$.

\end{itemize}
This shows that $\mathcal{E}(S_1, \Lambda_1) \cup \mathcal{E}(S_2, \Lambda_2)$ is a Riesz basis for $L^2[0,1)$ in all cases. 

\medskip

Finally, we assume that all intervals $I_1, I_2, I_3$ are {\it nonempty}.
Similarly as above, one can easily see that $\mathcal{E} (S_1 \backslash I_3 , \Lambda_1) \cup \mathcal{E} (S_2 \backslash I_3, \Lambda_2)$ is a Riesz basis for $L^2(I_1 \cup I_2)$; note that $\mathcal{E} (I_1 \cup I_2, \Lambda_1 \cup \Lambda_2)$ is a Riesz basis for $L^2(I_1 \cup I_2)$ by \cite[Proposition 2.1]{MM09} or \cite[Proposition 5.4]{BCMS19}.

First, we consider some special cases.

\begin{itemize}
\item[Case ({\romannumeral 1}).]
There exists some $k \in \{ 1, 2, 3 \}$ such that $\mathcal{L}_k = \{ k \}$.

By symmetry, we may assume that $k = 3$, that is, $\mathcal{L}_3 = \{ 3 \}$.
This means that $S_3 = I_3$.

\medskip

\noindent
{\it Injectivity of $T$.} \
Let $\bc = \{ c_n \}_{n \in \Z} \in \ell_2(\Z)$ be such that 
for a.e.~$x \in [0,1)$, 
\[
\sum_{\lambda\in \Lambda_1} c_\lambda \, e^{2\pi i \lambda x} \chi_{S_1} (x) + \sum_{\lambda\in \Lambda_2} c_\lambda \, e^{2\pi i \lambda x} \chi_{S_2} (x) + \sum_{\lambda\in \Lambda_3} c_\lambda \, e^{2\pi i \lambda x} \chi_{I_3} (x) = 0 . 
\]
Note that $\sum_{\lambda\in \Lambda_1} c_\lambda \, e^{2\pi i \lambda x} \chi_{S_1} (x) + \sum_{\lambda\in \Lambda_2} c_\lambda \, e^{2\pi i \lambda x} \chi_{S_2} (x) = 0$ for a.e.~$x \in I_1 \cup I_2$.
Since $\mathcal{E} (S_1 \backslash I_3 , \Lambda_1) \cup \mathcal{E} (S_2 \backslash I_3, \Lambda_2)$ is a Riesz basis for $L^2(I_1 \cup I_2)$, we obtain $c_\lambda = 0$ for $\lambda \in \Lambda_1 \cup \Lambda_2$.
In turn, since $\mathcal{E} (I_3, \Lambda_3)$ is a Riesz basis for $L^2(I_3)$, we get $c_\lambda = 0$ for $\lambda \in \Lambda_3$.

\medskip

\noindent
{\it Surjectivity of $T$.} \
Fix any function $f \in L^2[0,1)$.
Since $\mathcal{E} (S_1 \backslash I_3, \Lambda_1) \cup \mathcal{E} (S_2 \backslash I_3, \Lambda_2)$ is a Riesz basis for $L^2(I_1 \cup I_2)$, there exists a unique sequence $\{ c_\lambda \}_{\lambda \in \Lambda_1 \cup \Lambda_2} \in \ell_2(\Lambda_1 \cup \Lambda_2)$ such that
$\sum_{\lambda\in \Lambda_1} c_\lambda \, e^{2\pi i \lambda x} \chi_{S_1} (x) + \sum_{\lambda\in \Lambda_2} c_\lambda \, e^{2\pi i \lambda x} \chi_{S_2} (x) = f(x)$ for a.e.~$x \in I_1 \cup I_2$.
Define $f_3 \in L^2(I_3)$ by $f_3(x) := f(x) -  \sum_{\lambda\in \Lambda_1} c_\lambda \, e^{2\pi i \lambda x} \chi_{S_1} (x) - \sum_{\lambda\in \Lambda_2} c_\lambda \, e^{2\pi i \lambda x} \chi_{S_2} (x)$ for $x \in I_3$.
Since $\mathcal{E} (I_3, \Lambda_3)$ is a Riesz basis for $L^2(I_3)$, there exists a unique sequence $\{ c_\lambda \}_{\lambda \in \Lambda_3} \in \ell_2(\Lambda_3)$ such that $\sum_{\lambda\in \Lambda_3} c_\lambda \, e^{2\pi i \lambda x} (x) = f_3(x)$ for a.e.~$x \in I_3$.
Then it holds for a.e.~$x \in [0,1)$, 
\[
\sum_{\lambda\in \Lambda_1} c_\lambda \, e^{2\pi i \lambda x} \chi_{S_1} (x) + \sum_{\lambda\in \Lambda_2} c_\lambda \, e^{2\pi i \lambda x} \chi_{S_2} (x) + \sum_{\lambda\in \Lambda_3} c_\lambda \, e^{2\pi i \lambda x} \chi_{I_3} (x)  = f(x)  . 
\]
Hence, the system $\mathcal{E} (S_1, \Lambda_1) \cup \mathcal{E} (S_2, \Lambda_2) \cup \mathcal{E} (S_3, \Lambda_3)$ is a Riesz basis for $L^2[0,1)$.

\item[Case ({\romannumeral 2}).]
There exists some $k \in \{ 1, 2, 3 \}$ such that $k \notin \mathcal{L}_\ell$ for all $\ell \neq k$.

\medskip

By symmetry, we may assume that $k = 3$, that is, $3 \notin \mathcal{L}_1$ and $3 \notin \mathcal{L}_2$. This means that $I_3$ is contained in $S_3$ but not in $S_1$ and $S_2$.

\medskip

\noindent
{\it Injectivity of $T$.} \
Let $\bc = \{ c_n \}_{n \in \Z} \in \ell_2(\Z)$ be such that
for a.e.~$x \in [0,1)$, 
\[
\sum_{\lambda\in \Lambda_1} c_\lambda \, e^{2\pi i \lambda x} \chi_{S_1} (x) + \sum_{\lambda\in \Lambda_2} c_\lambda \, e^{2\pi i \lambda x} \chi_{S_2} (x) + \sum_{\lambda\in \Lambda_3} c_\lambda \, e^{2\pi i \lambda x} \chi_{S_3} (x) = 0 . 
\]
Note that for a.e.~$x \in I_3$, we have $\sum_{\lambda\in \Lambda_3} c_\lambda \, e^{2\pi i \lambda x} = 0$.
Since $\mathcal{E} (I_3, \Lambda_3)$ is a Riesz basis for $L^2(I_3)$, we obtain $c_\lambda = 0$ for $\lambda \in \Lambda_3$.
In turn, since $\mathcal{E} (S_1, \Lambda_1) \cup \mathcal{E} (S_2, \Lambda_2) = \mathcal{E} (S_1 \backslash I_3, \Lambda_1) \cup \mathcal{E} (S_2 \backslash I_3 , \Lambda_2)$ is a Riesz basis for $L^2(I_1 \cup I_2)$, we get
$c_\lambda = 0$ for $\lambda \in \Lambda_1 \cup \Lambda_2$.

\medskip

\noindent
{\it Surjectivity of $T$.} \
Likewise, this can be seen by first using that $\mathcal{E} (I_3, \Lambda_3)$ is a Riesz basis for $L^2(I_3)$, and then using that $\mathcal{E} (S_1, \Lambda_1) \cup \mathcal{E} (S_2, \Lambda_2)$ is a Riesz basis for $L^2(I_1 \cup I_2)$.

Hence, the system $\mathcal{E} (S_1, \Lambda_1) \cup \mathcal{E} (S_2, \Lambda_2) \cup \mathcal{E} (S_3, \Lambda_3)$ is a Riesz basis for $L^2[0,1)$.

\item[Case ($*$).]
The sets $(\mathcal{L}_k)^c = \{1,2,3\} \backslash \mathcal{L}_k$, $k=1,2,3$, are disjoint.

\medskip

For $\bc = \{ c_n \}_{n \in \Z} \in \ell_2(\Z)$, we compute
\begin{equation}\label{eqn:CaseStar}
\begin{split}
& \Big\| \sum_{n \in \Z} c_n \, e_n  - T(\bc) \Big\|_{L^2[0,1)}^2 \\
&= \Big\| \sum_{k=1}^3 \sum_{\lambda\in \Lambda_k} c_\lambda \, e^{2\pi i \lambda (\cdot)} \chi_{[0,1) \backslash S_k} \Big\|_{L^2[0,1)}^2 \\
&= \sum_{k=1}^3 \Big\| \sum_{\lambda\in \Lambda_k} c_\lambda \, e^{2\pi i \lambda (\cdot)} \Big\|_{L^2([0,1) \backslash S_k)}^2
\qquad \big(\text{since} \;\; [0,1) \backslash S_k \;\; \text{are disjoint}\big) \\
&\leq \sum_{k=1}^3 \Big( \Big\| \sum_{\lambda\in \Lambda_k} c_\lambda \, e^{2\pi i \lambda (\cdot)} \Big\|_{L^2[0,1)}^2 - \Big\| \sum_{\lambda\in \Lambda_k} c_\lambda \, e^{2\pi i \lambda (\cdot)} \Big\|_{L^2(I_k)}^2 \Big) \\
&\leq \sum_{k=1}^3 \Big( \sum_{\lambda\in \Lambda_k} |c_\lambda|^2 - \alpha_k \sum_{\lambda\in \Lambda_k} |c_\lambda|^2 \Big) \\
&\leq \big(1- \min\{ \alpha_1,\alpha_2,\alpha_3  \}\big) \, \| \bc \|_{\ell_2}^2 ,
\end{split}
\end{equation}
where $0 < \alpha_k \leq |I_k| < 1$ is a lower Riesz bound of $\mathcal{E} (I_k, \Lambda_k)$ in $L^2(I_k)$.
Using Theorem~\ref{thm:PaleyWienerStabilityGeneral} with $\lambda := (1- \min\{ \alpha_1,\alpha_2,\alpha_3  \})^{1/2} < 1$, we conclude that
$\mathcal{E} (S_1, \Lambda_1) \cup \mathcal{E} (S_2, \Lambda_2) \cup \mathcal{E} (S_3, \Lambda_3)$ is a Riesz basis for $L^2[0,1)$.
\end{itemize}

We now consider the full set of cases. There are $2^{2 \cdot 3} = 64$ possible choices for $(\mathcal{L}_1, \mathcal{L}_2, \mathcal{L}_3)$.

\begin{itemize}
\item[Case 1.]
$\mathcal{L}_3 = \{ 3 \}$.
This corresponds to Case ({\romannumeral 1}).
Note that there are 16 subcases; $\mathcal{L}_1$ is one of $\{ 1 \}$, $\{ 1,2 \}$, $\{ 1,3 \}$, $\{ 1,2,3 \}$, 
and $\mathcal{L}_2$ is one of $\{ 2 \}$, $\{ 1,2 \}$, $\{ 2,3 \}$, $\{ 1,2,3 \}$, independently.

\item[Case 2.]
$\mathcal{L}_3 = \{ 1, 3 \}$.

\begin{itemize}
\item[2-1.]
$\mathcal{L}_1 = \{ 1 \}$, $\mathcal{L}_2 = \{ 2 \}$, $\mathcal{L}_3 = \{ 1, 3 \}$.
This corresponds to Case ({\romannumeral 1}).

\item[2-2.]
$\mathcal{L}_1 = \{ 1 \}$, $\mathcal{L}_2 = \{ 1,2 \}$, $\mathcal{L}_3 = \{ 1, 3 \}$.
This corresponds to Case ({\romannumeral 1}).

\item[2-3.]
$\mathcal{L}_1 = \{ 1 \}$, $\mathcal{L}_2 = \{ 2,3 \}$, $\mathcal{L}_3 = \{ 1, 3 \}$.
This corresponds to Case ({\romannumeral 1}).

\item[2-4.]
$\mathcal{L}_1 = \{ 1 \}$, $\mathcal{L}_2 = \{ 1,2,3 \}$, $\mathcal{L}_3 = \{ 1, 3 \}$.
This corresponds to Case ({\romannumeral 1}).

\item[2-5.]
$\mathcal{L}_1 = \{ 1,2 \}$, $\mathcal{L}_2 = \{ 2 \}$, $\mathcal{L}_3 = \{ 1, 3 \}$.
This corresponds to Case ({\romannumeral 1}).

\item[2-6.]
$\mathcal{L}_1 = \{ 1,2 \}$, $\mathcal{L}_2 = \{ 1, 2 \}$, $\mathcal{L}_3 = \{ 1, 3 \}$.
This corresponds to Case ({\romannumeral 2}) with $k = 3$.

\item[2-7.]
$\mathcal{L}_1 = \{ 1,2 \}$, $\mathcal{L}_2 = \{ 2, 3 \}$, $\mathcal{L}_3 = \{ 1, 3 \}$.
This corresponds to Case ($*$).

\item[2-8.]
$\mathcal{L}_1 = \{ 1,2 \}$, $\mathcal{L}_2 = \{ 1, 2, 3 \}$, $\mathcal{L}_3 = \{ 1, 3 \}$.
This corresponds to Case ($*$).

\item[2-9.]
$\mathcal{L}_1 = \{ 1,3 \}$, $\mathcal{L}_2 = \{ 2 \}$, $\mathcal{L}_3 = \{ 1, 3 \}$.
This corresponds to Case ({\romannumeral 1}).

\item[2-10.]
$\mathcal{L}_1 = \{ 1,3 \}$, $\mathcal{L}_2 = \{ 1, 2 \}$, $\mathcal{L}_3 = \{ 1, 3 \}$.
This corresponds to Case ({\romannumeral 2}) with $k = 2$.

\item[2-11.]
$\mathcal{L}_1 = \{ 1,3 \}$, $\mathcal{L}_2 = \{ 2, 3 \}$, $\mathcal{L}_3 = \{ 1, 3 \}$.
This corresponds to Case ({\romannumeral 2}) with $k = 2$.

\item[2-12.]
$\mathcal{L}_1 = \{ 1,3 \}$, $\mathcal{L}_2 = \{ 1, 2, 3 \}$, $\mathcal{L}_3 = \{ 1, 3 \}$.
This corresponds to Case ({\romannumeral 2}) with $k = 2$.

\item[2-13.]
$\mathcal{L}_1 = \{ 1,2,3 \}$, $\mathcal{L}_2 = \{ 2 \}$, $\mathcal{L}_3 = \{ 1, 3 \}$.
This corresponds to Case ({\romannumeral 1}).

\item[2-14.]
$\mathcal{L}_1 = \{ 1,2,3 \}$, $\mathcal{L}_2 = \{ 1, 2 \}$, $\mathcal{L}_3 = \{ 1, 3 \}$.
This corresponds to Case ($*$).

\item[2-15.]
$\mathcal{L}_1 = \{ 1,2,3 \}$, $\mathcal{L}_2 = \{ 2, 3 \}$, $\mathcal{L}_3 = \{ 1, 3 \}$.
This corresponds to Case ($*$).

\item[2-16.]
$\mathcal{L}_1 = \{ 1,2,3 \}$, $\mathcal{L}_2 = \{ 1, 2, 3 \}$, $\mathcal{L}_3 = \{ 1, 3 \}$.
This corresponds to Case ($*$).
\end{itemize}

\item[Case 3.]
$\mathcal{L}_3 = \{ 2, 3 \}$.

There are 16 subcases, each of which is symmetric with a subcase of Case 2.

\item[Case 4.]
$\mathcal{L}_3 = \{ 1, 2, 3 \}$.

There are 16 subcases.
The subcase $\mathcal{L}_1 = \mathcal{L}_2 = \mathcal{L}_3 = \{ 1,2,3 \}$ corresponds to Case ($*$).
In fact, this subcase is trivial, since $\mathcal{E} ([0,1) , \Z)$ is an orthonormal basis for $L^2[0,1)$.
All other subcases have a symmetric subcase in Cases 1, 2, 3.
\end{itemize}
\end{proof}

%

\section*{Acknowledgments} 
D.G.~Lee acknowledges support by the DFG (German Research Foundation) Grant PF 450/9-2.
G.E.~Pfander acknowledges support by the DFG Grant PF 450/11-1.

\end{document}